\documentclass[11pt]{amsart}
\usepackage{graphicx}
\usepackage{amscd}
\usepackage{amsmath}
\usepackage{psfig}
\usepackage{amsfonts}
\usepackage{amssymb}
\newtheorem{theorem}{Theorem}[section]
\theoremstyle{plain}

\newtheorem{corollary}{Corollary}[section]

\newtheorem{lemma}{Lemma}[section]

\newtheorem{proposition}{Proposition}[section]

\numberwithin{equation}{section}

\begin{document}
\title[Isometric Embedding]{Counterexamples for Local Isometric Embedding}
\author{Nikolai NADIRASHVILI}
\address{Department of Mathematics\\
University of Chicago\\
5734 S. University Ave., Chicago, IL 60637}
\email{nicholas@math.uchicago.edu}
\author{Yu YUAN}
\address{Department of Mathematics\\
University of Chicago\\
5734 S. University Ave., Chicago, IL 60637\\
and University of Washington\\
Seattle, WA 98195}
\email{yuan@math.uchicago.edu, yuan@math.washington.edu}
\thanks{Both authors are partially supported by NSF grants, and the second author also
by a Sloan Research Fellowship}
\date{July 3, 2002}
\maketitle

\section{ Introduction}

In this paper, we construct metrics on 2-manifold which cannot be even locally
isometrically embedded in the Euclidean space $\mathbb{R}^{3}.$ By isometric
embedding of $\left(  M^{2},g\right)  $ with $g=\sum_{i,j=1}^{2}g_{ij}%
dx_{i}dx_{j}$ in $\mathbb{R}^{3},$ we mean there exists a surface in
$\mathbb{R}^{3}$ with the induced metric equaling $g,$ namely, the three
coordinate functions $\left(  X\left(  x_{1},x_{2}\right)  ,Y\left(
x_{1},x_{2}\right)  ,Z\left(  x_{1},x_{2}\right)  \right)  $ defined on
$M^{2}$ satisfy
\[
dX^{2}+dY^{2}+dZ^{2}=\sum_{i,j=1}^{2}g_{ij}dx_{i}dx_{j}.
\]

To be precise, we state the results in the following

\begin{theorem}
There exists a smooth metric $g$ in $B_{1}\subset\mathbb{R}^{2}$ with Gaussian
curvature $K_{g}\leq0$ such that there is no $C^{3}$ isometric embedding of
$\left(  B_{r}\left(  0\right)  ,g\right)  $ in $\mathbb{R}^{3}$ for any $r>0.$
\end{theorem}

\begin{theorem}
There exists a smooth metric $g$ in $B_{1}\subset\mathbb{R}^{2}$ with Gaussian
curvature $K_{g}(0)=0$ and $K_{g}\left(  x\right)  <0$ for $x\neq0$ such that
there is no $C^{3,\alpha}$ isometric embedding of $\left(  B_{r}\left(
0\right)  ,g\right)  $ in $\mathbb{R}^{3}$ for any $r>0$ and $\alpha>0.$
\end{theorem}

Pogorelov [P2] constructed a simple $C^{2,1}$ metric $g$ in $B_{1}\subset
\mathbb{R}^{2}$ with sign-changing Gaussian curvature such that $\left(  
B_{r}%
,g\right)  $ cannot be realized as a $C^{2}$ surface in $\mathbb{R}^{3}$ for
any $r>0.$ Recently the first author [N] gave a $C^{\infty}$ metric $g$ on
$B_{1}$ with no smooth isometric embedding of $\left(  B_{r},g\right)  $ in
$\mathbb{R}^{3}$ for any $r>0.$ The sign of the Gaussian curvature $K_{g}$
also changes.

On the positive side, when the sign of $K_{g}$ for any smooth metric $g$ does
not change, the local smooth isometric embedding was settled by Pogorelov
[P1], Nirenberg [Ni], and Hartman and Winter [HW2]. When $K_{g}\geq0$ for the
$C^{k}$ metric with $k\geq10,$ there is a $C^{k-6}$ isometric embedding of
$\left(  B_{r_{k}},g\right)  $ in $\mathbb{R}^{3},$ this was done by Lin [L1].
When $K_{g}$ changes sign cleanly, namely, $K_{g}\left(  0\right)  =0,\nabla
g\left(  0\right)  \neq0$ for a $C^{k}$ metric $g,$ Lin [L2] showed that there
exists a $C^{k-3}$ isometric embedding in $\mathbb{R}^{3}$ for $\left(
B_{r_{k}},g\right)  $ with $k\geq6.$ When $K_{g}\leq0$ and $\nabla^{2}%
K_{g}\left(  0\right)  \neq0$ for the smooth metric $g,$ there is a local
smooth isometric embedding of $g$ in $\mathbb{R}^{3},$ see Iwasaki [I]. When
$K_{g}=-{x_{1}^{2m}}\widetilde{K}(x)$ with $\widetilde{K}(0)>0$ for the smooth
metric $g$, the same local isometric embedding also holds, see Hong [H].
Recently, Han, Hong, and Lin [HHL] showed that the local isometric embedding
exists under the assumption $K_{g}\leq0$ with a certain non-degeneracy of the
gradient of $K_{g},$ or $K_{g}\leq0$ with finite order vanishing.

If one allows higher dimensional ambient space, say $\mathbb{R}^{4}$, Poznyak
[Po1] proved that any smooth metric $g$ on $M^{2}$ can be locally smoothly
isometrically embedded in $\mathbb{R}^{4}.$ In fact, any $C^{k}$ metric on
n-manifold $M^{n}$ has a $C^{k}$ global isometric embedding in $\mathbb{R}%
^{N_{n}}$ with $N_{n}$ large for $3\leq k\leq\infty.$ This is the work by Nash [Na2].

If we start with an analytic metric $g$ on $M^{n}$, one always has a local
analytic isometric embedding of $\left(  M^{n},g\right)  $ in $\mathbb{R}%
^{n\left(  n+1\right)  /2}.$ This was proved by Janet [J], Cartan [C] very
earlier on, and initiated by Schlaefli in 1873!

Lastly, any $C^{0}$ metric $g$ on a compact n-manifold $M^{n}$ which can be
differentially embedded in $\mathbb{R}^{n+1}$ has a $C^{1}$ isometric
embedding in $\mathbb{R}^{n+1},$ see Nash [Na1] and Kuiper [K].

For general description and further results on isometric embedding problem, we
refer to [GR], [P2] and [Y].

The heuristic idea of the construction is to arrange the metric $g$ in $B_{1}$
so that the second fundamental form of any isometric embedded surface in
$\mathbb{R}^{3}$, II$\circ i$ vanishes at one point, where $i:\left(
B_{1},g\right)  \rightarrow \mathbb{R}^{3}$ is the isometric embedding 
which is
supposed to exist. Further we force II$\circ i$ to vanish along the boundary
of a small domain $\Omega$ near the center of $B_{1},$ where the Gaussian
curvature $K_{g}<0$ (in $\Omega$). By the maximal principle, one cannot have a
saddle surface with vanishing second fundamental form along the boundary. So
$\left(  \Omega,g\right)  $ cannot be realized in $\mathbb{R}^{3}.$ We repeat
the construction near the center of $B_{1}$ at every scale so that $\left(
B_{1},g\right)  $ is not isometrically embeddable in $\mathbb{R}^{3}$ near the center.

The way to force II$\circ i$ to vanish at one point, say o, is the following.
We modify the flat metric $g_{0}=dx^{2}$ in $\mathbb{R}^{2}$ only over certain
region $\Lambda$ slightly away from the center o to a new one $g$ so that, for
a segment $A_{1}A_{2}$ with $A_{1},$ $A_{2}\in\partial\Lambda,$ the length of
$A_{1}A_{2}$ under $g$ is shorter than the one of the geodesic $A_{1}A_{2}$
under the flat $g_{0},$ and $K_{g}\leq0$ in a subregion $\Lambda_{s}$
containing $A_{1}A_{2}$. Because of $\det$II$\left(  i\left(  0\right)
\right)  =0,$ we only need to deal with the other principle curvature. Suppose
the second one $\kappa_{2}\neq0,$ say $\kappa_{2}<0.$ We show that there is a
flat concave cylinder $\Sigma$ near $i\left(  B_{1}\right)  ,$ which is
isometric to $\left(  B_{1},g_{0}\right)  $ provided the embedding $i$ is
$C^{3}$ (This assertion for $C^{2}$ embedding case remains unclear to us). 
Now $i\left(  A_{1}A_{2}\right)  $ supported on the saddle surface
$i\left(  \Lambda_{s}\right)  $ can only stay above the concave cylinder
$\Sigma.$ Then the length of $i\left(  A_{1}A_{2}\right)  $ is longer than the
one of the projection of $i\left(  A_{1}A_{2}\right)  $ down to the flat
$\Sigma,$ call it $P\circ i\left(  A_{1}A_{2}\right)  .$ We know the length of
$P\circ i\left(  A_{1}A_{2}\right)  $ under $g_{0}$ is equal to or longer than 
that of
the geodesic $A_{1}A_{2}$ under $g_{0}.$ But we start from $A_{1}A_{2}$ with
shorter length under $g$ than under $g_{0}.$ This contradiction shows that
II$\circ i\left(  0\right)  $ vanishes.

Inevitably, $K_{g}$ is positive somewhere in $\Lambda$ if $\Lambda$ is
surrounded by flat region with metric $dx^{2}.$ We add ``tails'' extending to
the boundary $\partial B_{1}$ for the modifying regions $\Lambda,$ modify the
metric on the tails, then we have the $g$ with $K_{g}\leq0$ in $B_{1}.$ It
turns out that we cannot work with a segment in the construction, we go with a
minimal tree connecting three points on $\partial\Lambda$ for each $\Lambda,$
see section 2 for details.

Now that we have a non-isometrically embeddable metric (with nonpositive
Gaussian curvature), the nearby metrics are almost non-isometrically
embeddable. Based on this observation, we construct a non-isometrically
embeddable metric with negative Gaussian curvature except for one point in
section 3.

\section{Metric with nonpositive curvature}

Recall any three segments in $\mathbb{R}^{2}$ with equal angles $\frac{2}%
{3}\pi$ at the common vertex form a minimal tree $T,$ namely, the length of
$T$ is less than that of any arcs connecting the other three vertices.

\begin{lemma}
\label{Moon} Let $u=-\operatorname{Im}e^{\log^{2}z}=-e^{\log^{2}r-\theta^{2}%
}\sin\left(  2\theta\log r\right)  ,$ $0<\theta<2\pi.$ Then there exists a
large integer $K$ such that%
\[
\int_{T}uds<0,
\]
where the minimal tree $T=AA_{1}\cup AA_{2}\cup AA_{3}$ with $A=\left(
-e^{-K},0\right)  ,$ $A_{2}=\left(  -1,0\right)  ,$ $A_{1,}$ $A_{2}\in\partial
B_{1},$ $\angle A_{1}AA_{2}=\angle A_{2}AA_{3}=\frac{2}{3}\pi.$ Moreover,
$u_{r}<0$ for $r=1.$
\end{lemma}

\begin{proof}
Set $\Omega_{u}=B_{1}\cap$Sector$A_{1}AA_{2},$ $\Omega_{l}=B_{1}\cap
$Sector$A_{2}AA_{3},$ $\widehat{A_{1}A_{2}}=\partial\Omega_{u}\cap\partial
B_{1},$ $\widehat{A_{2}A_{3}}=\partial\Omega_{l}\cap\partial B_{1}.$ Let the
angle from $A_{1}A$ to $x$ be $\varphi,$ or $\varphi\left(  x\right)  =\angle
A_{1}Ax,$ then $0\leq\varphi\left(  x\right)  \leq\frac{4}{3}\pi$ for
$x\in\Omega_{u}\cup\Omega_{l}.$

We apply Green formula to harmonic functions $u$ and $\varphi$ in $\Omega_{u}$
and $\Omega_{l},$%
\begin{align*}
\int_{\partial\Omega_{u}}u\varphi_{\gamma}ds &  =\int_{\partial\Omega_{u}%
}\varphi u_{\gamma}ds\\
\int_{\partial\Omega_{l}}u\left(  \varphi-\frac{4}{3}\pi\right)  _{\gamma}ds
&  =\int_{\partial\Omega_{l}}\left(  \varphi-\frac{4}{3}\pi\right)  u_{\gamma
}ds,
\end{align*}
where $\gamma$ is the outward unit normal of the integral domain. We then have%
\begin{align*}
\int_{AA_{1}}-uds+\int_{AA_{2}}uds &  =\int_{\widehat{A_{1}A_{2}}}\varphi
u_{r}ds+\int_{AA_{2}}\frac{2}{3}\pi u_{\theta}ds\\
\int_{AA_{2}}-uds+\int_{AA_{3}}uds &  =\int_{\widehat{A_{2}A_{3}}}\left(
\varphi-\frac{4}{3}\pi\right)  u_{r}ds+\int_{AA_{2}}\frac{2}{3}\pi u_{\theta
}ds.
\end{align*}
It follows that%
\begin{align*}
\int_{AA_{1}\cup AA_{3}}uds &  =2\int_{AA_{2}}uds+\int_{\widehat{A_{1}A_{2}}%
}-\varphi u_{r}ds+\int_{\widehat{A_{2}A_{3}}}\left(  \varphi-\frac{4}{3}%
\pi\right)  u_{r}ds\\
&  =2\int_{AA_{2}}uds+\int_{\widehat{A_{1}A_{2}}}\varphi e^{-\theta^{2}%
}2\theta ds+\int_{\widehat{A_{2}A_{3}}}\left(  \frac{4}{3}\pi-\varphi\right)
e^{-\theta^{2}}2\theta ds.
\end{align*}
On the other hand,%
\begin{align*}
\int_{AA_{2}}uds &  =\int_{e^{-K}}^{e^{0}}-e^{\left(  \log^{2}r-\pi
^{2}\right)  }\sin\left(  2\pi\log r\right)  dr\\
&  =\frac{1}{2\pi e^{\pi^{2}}}\int_{-2\pi K}^{0}-e^{\left(  \frac{t^{2}}%
{4\pi^{2}}+\frac{t}{2\pi}\right)  }\sin tdt.
\end{align*}
We choose large enough integer $K$ so that $\int_{AA_{2}}uds<0$ and
\[
2\int_{AA_{2}}uds+\int_{\widehat{A_{1}A_{2}}}\varphi e^{-\theta^{2}}2\theta
ds+\int_{\widehat{A_{2}A_{3}}}\left(  \frac{4}{3}\pi-\varphi\right)
e^{-\theta^{2}}2\theta ds<0.
\]
Therefore%
\[
\int_{T}uds<0.
\]
\end{proof}

\noindent\textbf{Remark.} By applying Green formula to the above harmonic
function $u$ and linear functions, one sees that $\int_{\Gamma}uds>0$ for any
segment $\Gamma\subset\Omega_{u}\cup\Omega_{l},$ connecting two boundary
points on $\partial B_{1}.$

\begin{lemma}
\label{Tail}There exists a function $v\in C_{0}^{\infty}\left(  B_{1.1}%
\right)  $ satisfying%
\begin{align*}
v &  =0\;\;\;\;\;\;in\;\;\;\;\left\{  \left(  x_{1},x_{2}\right)
|x_{1}<0.9\right\}  \backslash B_{1}\\
\bigtriangleup v &  \geq0\;\text{\ \ \ \ \ \ }\;in\;\;\;\;\;\;\;\;\;B_{1}\\
\int_{T}vds &  <0
\end{align*}
where the minimal tree $T=CC_{1}\cup CA_{2}\cup CC_{3}$ with $A_{2}=\left(
-1,0\right)  ,$ $C=\left(  -\frac{1}{10}e^{-K}-0.8,0\right)  ,$ $C_{1},$
$C_{3}\in\partial B_{1}$ and $\angle C_{1}CA_{2}=\angle A_{2}CC_{3}=\frac
{2}{3}\pi.$ Moreover $T\subset\left\{  \left(  x_{1},x_{2}\right)
|x_{1}<-0.1\right\}  .$
\end{lemma}%

\begin{figure}
[ptbh]
\begin{center}
\includegraphics[
trim=0.000000in 0.000000in -0.002490in 0.000886in,
height=6.4954cm,
width=6.0934cm
]%
{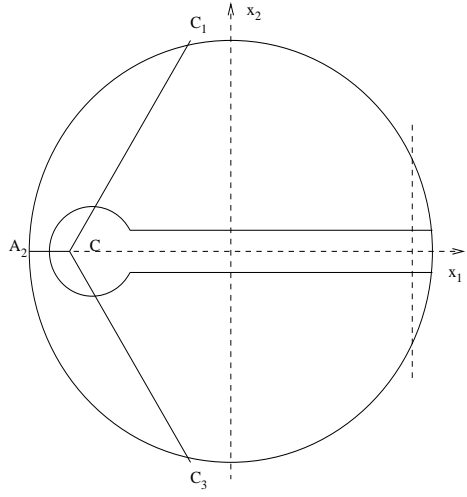}%
\caption{Minimal tree inside the half ball.}%
\end{center}
\end{figure}

\begin{proof}
Set \ $D=\left(  -e^{-2K},0\right)  ,$ $D_{1},$ $D_{2}\in\partial B_{1}$ with
$\angle D_{1}DA_{2}=\angle A_{2}DD_{3}=\frac{2}{3}\pi,$ and $D_{4}=\left(
20,x_{2}\left(  D_{3}\right)  \right)  ,$ $D_{5}=\left(  20,x_{2}\left(
D_{1}\right)  \right)  .$ Set $\Omega_{p}=$Pentagon$D_{1}DD_{3}D_{4}D_{5}.$
Let $w$ satisfy%
\begin{align*}
\bigtriangleup w &  =0\;\;\;\;in\;\;\;\Omega_{p}\;\;\\
w &  =u\;\;\;\;on\;\;\;D_{1}D\cup D_{3}D\\
w &  =0\;\;\;\;on\;\;\;D_{1}D_{5}\cup D_{3}D_{4}\\
w &  =N\;\;\;\;on\;\;\;D_{4}D_{5}\\
w &  =u\;\;\;\;\;in\;\;\;\;B_{1}\backslash\text{Sector}D_{1}DD_{3},
\end{align*}
where $u$ is the one in Lemma \ref{Moon}.

We choose large enough $N$ so that $w_{\gamma}>u_{\gamma}$ on $D_{1}D\cup
D_{3}D$ and $w_{\gamma}>0$ on $D_{1}D_{5}\cup D_{3}D_{4},$ where $\gamma$ is
the inward unit normal of $\partial\Omega_{p}$ this time. (If one insists, we
can smooth off $\partial\Omega_{p}.$)

Next we mollify $w$ by the usual (radially symmetric) mollifier $\rho_{\delta
}\in C_{0}^{\infty}\left(  B_{\delta}\right)  $ with $0<\delta<e^{-2K}$ to be
determined later. We see that the smooth function $w\ast\rho_{\delta}$
satisfies%
\begin{align*}
\bigtriangleup w\ast\rho_{\delta}\left(  x\right)   &  \geq
0\;\;\;\;\;for\;\;\;x_{1}\leq19.9\\
w\ast\rho_{\delta}\left(  x\right)   &  =u\;\;\;\;\;for\;\;\;x\;\text{\ inside
}\Omega_{i}=B_{1}\backslash\text{Sector}D_{1}DD_{3}\;\text{and }%
\delta\;\text{away from }\partial\Omega_{i}\\
w\ast\rho_{\delta}\left(  x\right)   &
=0\;\;\;\;\;for\;\;\;x\;\text{\ outside }\Omega_{o}=\left(  B_{1}%
\backslash\text{Sector}D_{1}DD_{3}\right)  \cup\Omega_{p}\;\text{and }%
\delta\;\text{away from }\partial\Omega_{o}.
\end{align*}
Finally, set $C_{0}=\left(  -0.8,0\right)  $ and%
\[
v\left(  x\right)  =w\ast\rho_{\delta}\left(  10\left(  x-C_{0}\right)
\right)  .
\]
By making $\delta$ even smaller yet positive if necessary so that $\int
_{T}vds<0$, we obtain the desired function $v$ in the above lemma.
\end{proof}

\begin{corollary}
\label{Minimal tree}Let $v$ be the function in Lemma \ref{Tail}. There exists
a family of smooth metrics in $\mathbb{R}^{2}$%
\[
g_{\delta}=e^{2\delta v}dx^{2}\;\;\;\;\text{for \ }0<\delta<\delta_{0}%
\]
such that%
\begin{align*}
g_{\delta} &  =dx^{2}\;\;\;\;\;\text{in \ \ }\left\{  \left(  x_{1}%
,x_{2}\right)  \left|  x_{1}<0.9\right|  \right\}  \backslash B_{1}\\
K_{g_{\delta}} &  \leq0\;\;\;\;\;\;\;\;\text{in \ \ }B_{1}\\
L\left(  T,g_{\delta}\right)   &  <L\left(  T,dx^{2}\right)  ,
\end{align*}
where $L\left(  T,g\right)  $ is the length of the minimal tree $T$ from Lemma
\ref{Tail} in metric $g.$
\end{corollary}

\begin{proof}
We only prove the last two inequalities. One has%
\[
K_{g_{\delta}}=-e^{-2\delta v}\bigtriangleup\left(  \delta v\right)
\leq0\;\;\;\;\text{in \ \ \ }B_{1}.
\]
Also%
\begin{align*}
L\left(  T,g_{\delta}\right)   &  =\int_{T}e^{\delta v}ds\\
\left.  \frac{dL}{d\delta}\right|  _{\delta=0} &  =\int_{T}vds<0.
\end{align*}
Thus there exists $\delta_{0}$ such that $L\left(  T,g_{\delta}\right)
<L\left(  T,dx^{2}\right)  \;\;\;\;$for \ $0<\delta<\delta_{0}.$
\end{proof}

Let $\psi\in C^{1}\left(  \left[  -1,1\right]  \right)  $ satisfy $0\leq
\psi\leq1$ and $\psi\left(  \pm1\right)  =0.$ Set%
\begin{align*}
\gamma &  =\left\{  \left(  x_{1},x_{2}\right)  |x_{1}=\psi\left(
x_{2}\right)  ,\;\left|  x_{2}\right|  \leq1\right\}  ,\;Q=\left\{  \left(
x_{1},x_{2}\right)  |0<x_{1}<\psi\left(  x_{2}\right)  ,\;\left|
x_{2}\right|  \leq1\right\} \\
\Pi &  =\left[  0,2\right]  \times\left[  -2,2\right]  \subset R^{2}%
,\;\;\;\;\;\;\;\;\;\;\;\;\;\;\;\;F=\Pi\backslash Q.
\end{align*}

\begin{lemma}
\label{Extension}Let $f\in C^{3}\left(  F\right)  .$ Assume the graph $\Sigma$
of $f$ is flat or $\det D^{2}f$ $=0$ and $D^{2}f\neq0$ in $F.$ Also assume a
unit $C^{1}$ continuous eigenvector $V_{0}$ for the zero eigenvalue of
$D^{2}f$ is transversal to $\gamma.$ For any $0<\tau<1,$ there exists
$\varepsilon>0$ so that if $\left\|  D^{2}f-\left[
\begin{array}
[c]{cc}%
0 & 0\\
0 & -\tau
\end{array}
\right]  \right\|  \leq\varepsilon\tau,$ one can extend $f$ to $\Pi$ with the
graph of the extension being flat and concave.
\end{lemma}

\begin{proof}
We take the $C^{2}$ Legendre coordinate system on $F\subset\Pi$ (cf. [HW1]).%
\[
\left\{
\begin{array}
[c]{l}%
t=x_{1}\\
s=f_{2}\left(  x_{1},x_{2}\right)  .
\end{array}
\right.
\]
Notice that the graph of $f$, $\Sigma$ is flat, or $\det D^{2}f$ $=0,$ it
follows that $\left\{  \left(  x_{1},x_{2}\right)  |f_{2}\left(  x_{1}%
,x_{2}\right)  =s=const\right\}  $ is a straight segment in $\mathbb{R}^{2}$
and $x_{t}\left(  t,s\right)  $ ($\Vert V_{0}$) is independent of $t.$ Also
$\frac{\partial f}{\partial t}\left(  x\left(  t,s\right)  \right)  $ is
independent of $t.$ Hence we can represent a portion $\Sigma^{p}$ of the graph
$\Sigma$ in the ruling form%
\[
\left(  x_{1},x_{2},x_{3}\right)  \left(  t,s\right)  =h\left(  t,s\right)
=c\left(  s\right)  +t\delta\left(  s\right)  =\left(  t,x_{2}\left(
t,s\right)  ,f\left(  t,x_{2}\left(  t,s\right)  \right)  \right)  ,
\]
where $c\left(  s\right)  ,\;\delta\left(  s\right)  \in C^{2}$ and $s\in
S=\left[  f_{2}\left(  2,2\right)  ,f_{2}\left(  2,-2\right)  \right]
,\;\;t\leq2.$

We may assume $\nabla f\left(  2,0\right)  =0.$ If $\varepsilon$ is chosen
small enough, then $\delta\left(  s\right)  \;\left(  \Vert V_{0}\right)  $ is
close to $\left(  1,0,0\right)  $ in $C^{1}$ norm. Take $\varepsilon$ small,
then
\[
\left\{  \left(  x_{1},x_{2},f\left(  x_{1},x_{2}\right)  \right)  |\left(
\left(  x_{1},x_{2}\right)  \in\gamma\right)  \right\}  \subset\partial
\Sigma^{p}.
\]
Set $U=\left\{  \left(  t,s\right)  |-1\leq t\leq2,\;s\in S\right\}  .$ Take
$\varepsilon$ small so that $\left\|  \delta\left(  s\right)  -\left(
1,0,0\right)  \right\|  _{C^{1}}$ small, then $\left(  t,s\right)  \in U$ is a
$C^{2}$ coordinate system for $\Pi.$

Now $\Sigma^{e}=h\left(  U\right)  $ is a $C^{2}$, flat, concave graph over a
domain $\Omega$ in $\mathbb{R}^{2}$ with $\Pi\subset\Omega.$ Indeed, the
normal of $\Sigma^{e}$ is%
\[
N=\frac{h_{t}\times h_{s}}{\left\|  h_{t}\times h_{s}\right\|  }.
\]
We know%
\begin{align*}
h_{t} &  =\left(  1,\frac{-f_{21}}{f_{22}},f_{1}+f_{2}\frac{-f_{21}}{f_{22}%
}\right)  \overset{\varepsilon\rightarrow0}{\longrightarrow}\left(
1,0,0\right)  \\
h_{s} &  =\left(  0,\frac{1}{f_{22}},\frac{f_{2}}{f_{22}}\right)
\overset{\varepsilon\rightarrow0}{\longrightarrow}\left(  0,\frac{-1}{\tau
},\frac{-s}{\tau}\right)  ,
\end{align*}
then $h_{t}\times h_{s}\overset{\varepsilon\rightarrow0}{\longrightarrow
}\left(  0,\frac{s}{\tau},\frac{-1}{\tau}\right)  .$ So $\Sigma^{e}$ is a
$C^{2}$ graph if we choose $\varepsilon$ small enough.

Next, the second fundamental form of $\Sigma^{e}$ is%
\begin{align*}
II &  =\left[
\begin{array}
[c]{cc}%
\left\langle h_{tt},N\right\rangle  & \left\langle h_{ts},N\right\rangle \\
\left\langle h_{st},N\right\rangle  & \left\langle h_{ss},N\right\rangle
\end{array}
\right]  \\
&  =\frac{1}{\left\|  h_{t}\times h_{s}\right\|  }\left[
\begin{array}
[c]{cc}%
0 & 0\\
0 & \left\langle c^{\prime\prime}+t\delta^{\prime\prime},\delta\times\left(
c^{\prime}+t\delta^{\prime}\right)  \right\rangle
\end{array}
\right]
\end{align*}
and the Gaussian curvature%
\[
K_{g}=0.
\]
Finally, the nonzero principle curvature of $\Sigma^{e}$%
\[
\kappa=\left[  \frac{\tau^{3}}{\left(  1+s^{2}\right)  ^{3/2}}+o\left(
\varepsilon\right)  \right]  \left\langle c^{\prime\prime}+t\delta
^{\prime\prime},\delta\times\left(  c^{\prime}+t\delta^{\prime}\right)
\right\rangle .
\]
On the other hand, from the graph representation of $\Sigma^{p}$,
$\kappa\overset{\varepsilon\rightarrow0}{\longrightarrow}-\tau/\left(
1+s^{2}\right)  ^{3/2}.$ So for $t$ in a certain range close to $2,$ say
$t\in\left[  1,2\right]  ,$ the quadratic function in terms of $t,$%
\[
\left\langle c^{\prime\prime}+t\delta^{\prime\prime},\delta\times\left(
c^{\prime}+t\delta^{\prime}\right)  \right\rangle =a_{0}+a_{1}t+a_{2}t^{2}%
\]
is close to $-1/\tau^{2}$ as $\varepsilon\rightarrow0$ . It follows that
$a_{0}+a_{1}t+a_{2}t^{2}$ is still close to $-1/\tau^{2}$ for $t\in\left[
-1,2\right]  ,$ if we choose $\varepsilon$ small enough. So $\Sigma^{e}$ is concave.
\end{proof}

\begin{lemma}
\label{Bernstein}Let $f$ be the extended function in Lemma \ref{Extension},
let $w\in C^{2}\left(  \Pi\right)  $ satisfy $w=f$ on $F,$ $\det D^{2}w\leq0$
in $\Pi,$ and $\left\|  D^{2}w-\left[
\begin{array}
[c]{cc}%
0 & 0\\
0 & -\tau
\end{array}
\right]  \right\|  $ $_{C^{1}}\leq\varepsilon\tau.$ Then%
\[
f\leq w\;\;in\;\;\;\Pi.
\]
\end{lemma}

\begin{proof}
Suppose there is a point $x^{\prime}=\left(  x_{1}^{\prime},x_{2}^{\prime
}\right)  \in M$ such that $w\left(  x^{\prime}\right)  <f\left(  x^{\prime
}\right)  $. We know $x_{2}^{\prime}\in\left(  -1,1\right)  .$ For simplicity,
we may assume%
\[
f\left(  x^{\prime}\right)  -w\left(  x^{\prime}\right)  =\sup_{x_{2}%
\in\left[  -1,1\right]  }\left[  f\left(  x_{1}^{\prime},x_{2}\right)
-w\left(  x_{1}^{\prime},x_{2}\right)  \right]  .
\]
Then $f_{2}\left(  x^{\prime}\right)  =w_{2}\left(  x^{\prime}\right)  .$ It
follows that the two tangent lines $l_{f}$, $l_{w}$ to $f$ and $w$ at
$x^{\prime}$ in the plane $\left\{  \left(  x_{1},x_{2},x_{3}\right)
|x_{1}=x_{1}^{\prime}\right\}  $ are parallel. Since $w\left(  x_{1}^{\prime
},\cdot\right)  $ is concave, $l_{w}$ is above $w.$

Let $T\subset\mathbb{R}^{3}$ be the tangent plane to the graph $\Sigma_{f}$ of
$f$ at $\left(  x^{\prime},f\left(  x^{\prime}\right)  \right)  .$ Let
$R=T\cap\Sigma_{f}.$ Then $R$ is a segment (ruling) transversal to $l_{f}.$
Let $\left(  x^{0},z^{0}\right)  \in R$ with $x^{0}\in F,$ then $z^{0}%
=f\left(  x^{0}\right)  =w\left(  x^{0}\right)  .$ Let $l_{0}\subset T$
through $\left(  x^{0},z^{0}\right)  $ with $l_{0}\Vert l_{w}$. By the
concavity of $f=w$ in $F,$ $l_{0}$ is above the graph $\Sigma_{w}$ of $w.$

Let $m\left(  x\right)  $ be the linear function with graph as the plane $E$
through $l_{w}$ and $l_{0}$. Let $V=\left\{  \left(  x_{1},x_{2}\right)
|x_{1}^{\prime}<x_{1}<2,\;\left|  x_{2}\right|  <2\right\}  .$ Because
$\Sigma_{w}$ is a ruling surface on $F$, then%
\[
w\left(  x\right)  \leq m\left(  x\right)  \;\;\;\text{on\ \ \ }\partial V.
\]
Note that $\det D^{2}w\leq0,$ by the maximum principle,%
\[
w\left(  x\right)  \leq m\left(  x\right)  \;\;\;\;\text{in}\;\;\;\;V.
\]
On the other hand, there is $\left(  x^{\ast},w\left(  x^{\ast}\right)
\right)  \in R$ with $x^{\ast}\in V$ such that%
\[
w\left(  x^{\ast}\right)  >m\left(  x^{\ast}\right)  .
\]
This contradiction completes the proof of the above lemma.
\end{proof}

Let $\textsl{r}$ be a rotation in $\mathbb{R}^{2}$ through an angle $1^{\circ
}$. Let $v$ be the function in Lemma \ref{Tail}, set%
\[
w\left(  x\right)  =\sum_{i=1}^{360}v\left(  \textsl{r}^{i}\left(
1000x\right)  -\left(  360,0\right)  \right)  .
\]
Pick two sequences $z_{n}\in\mathbb{R}^{2}$ and $\rho_{n}>0$ such that%
\[
z_{n}\longrightarrow0\;\;\;\text{as\ \ \ \ }n\longrightarrow+\infty
\]%
\[
B_{\rho_{n}}\left(  z_{n}\right)  \cap B_{\rho_{k}}\left(  z_{k}\right)
=\emptyset\;\;\;\text{for \ \ \ }n\neq k.
\]
Take another sequence $\delta_{n}>0$ going to $0$ fast enough so that the
smooth metric $g_{\text{II}}$ in $\mathbb{R}^{2}$ satisfying%
\begin{align*}
g_{\text{II}}  &  =e^{2\delta_{n}w\left(  z_{n}+x/\rho_{n}\right)  }%
dx^{2}\;\;\;\;\text{in \ \ }B_{\rho_{n}}\left(  z_{n}\right) \\
g_{\text{II}}  &  =dx^{2}%
\;\;\;\;\;\;\;\;\;\;\;\;\;\;\;\;\;\;\;\text{otherwise.}%
\end{align*}

\noindent\textbf{Remark.} Certainly our $v$ is only smooth in $B_{1.1}\left(
0\right)  $, that leaves the function $w$ nonsmooth, even undefined near the
corresponding tails. At this stage, we do not need any information on the
metric $g_{II}$ near those tails (Figure 1 and 3). We can make a smooth
extension of $v$ to $\mathbb{R}^{2}$ with $v\in C_{0}^{\infty}\left(
B_{2}\right)  $ if one insists. Then the Gaussian curvature of $g$ would be
positive near the transition region. In the proof of Theorem 1.1, we will
extend the tails to the boundary, make $v$ a smooth subharmonic function
inside the unit ball. Then the Gaussian curvature would be nonpositive in the
unit ball.

\begin{proposition}
\label{II=0}Let $i$ be a $C^{3}$ isometric embedding%
\[
i:\left(  B_{r}\left(  0\right)  ,g_{\text{II}}\right)  \longrightarrow
\mathbb{R}^{3}%
\]
for some $r>0.$ Then the second fundamental form of $i\left(  B_{r}\left(
0\right)  \right)  $ vanishes at $i\left(  0\right)  ,$ or \text{II}$\left(
i\left(  0\right)  \right)  =0.$
\end{proposition}

\begin{proof}
We may assume $i\left(  B_{r}\right)  $ is the graph $\Sigma_{w}$ of a
function $x_{3}=w\left(  x_{1},x_{2}\right)  $ and $w\left(  0\right)  =0$,
$\nabla w\left(  0\right)  =0.$ Then II$\left(  i\left(  0\right)  \right)
=D^{2}w\left(  0\right)  $ and $\det D^{2}w\left(  0\right)  =0.$ Suppose%
\[
D^{2}w\left(  0\right)  \neq0.
\]
Let $P_{3}$ be the projection from $\mathbb{R}^{3}$ to x$_{1}\_$x$_{2}$ plane.
Set $J\left(  x\right)  =P_{3}\left(  i\left(  x\right)  \right)  $. We may
assume DJ is the identity map on the tangent space $\mathbb{R}^{2}$ at $0,$
and%
\[
D^{2}w\left(  0,0\right)  =\left[
\begin{array}
[c]{cc}%
0 & 0\\
0 & -\tau
\end{array}
\right]  .
\]
For a sufficiently large $n,$\ $B_{\rho_{n}}\left(  z_{n}\right)  \subset
B_{r}$ and%
\[
g_{II}=e^{2\delta_{n}v\left(  r^{180}\left(  1000\left(  z_{n}+x/\rho
_{n}\right)  \right)  -\left(  360,0\right)  \right)  }dx^{2}%
\]
in the $179^{\circ}$ to $181^{\circ}$ section of the ball $B_{\rho_{n}}\left(
z_{n}\right)  .$

In order to simply the presentation, we work with the metric $g_{\delta_{n}%
}=e^{2\delta_{n}v\left(  x\right)  }dx^{2}$ as in the Corollary \ref{Minimal
tree}. Let $\Sigma^{e}$ be the flat, concave extension of $i\left(  B_{2}%
^{-}\backslash B_{1}^{-}\right)  $ by Lemma \ref{Extension}, where $B_{\rho
}^{-}=\left\{  \left(  x_{1},x_{2}\right)  |x_{1}<0\right\}  \cap B_{\rho}.$
Note that we may consider the graph $x_{3}=w_{\varepsilon}\left(  x\right)
=w\left(  \varepsilon x\right)  $ for small $\varepsilon$, then%
\[
\left\|  D^{2}w_{\varepsilon}-\left[
\begin{array}
[c]{cc}%
0 & 0\\
0 & -\varepsilon^{2}\tau
\end{array}
\right]  \right\|  _{C^{1}}\leq\varepsilon^{3},
\]
make the extension, then scale back.

Since $i\left(  B_{1}^{-}\right)  $ is negatively curved, or $\det D^{2}%
w\leq0$ and concave, we apply Lemma \ref{Bernstein} to conclude that $i\left(
B_{1}^{-}\right)  $ is above $\Sigma^{e}.$

Let $P$ be the normal projection of points $p$ above $\Sigma^{e}$ down to
$\Sigma^{e},$ that is $\left[  p-P\left(  p\right)  \right]  \perp\Sigma^{e}.$
By concavity of $\Sigma^{e}$, we have
\[
\text{Length}\left(  T,g_{\delta_{n}}\right)  =\text{Length}\left(  i\left(
T\right)  ,g_{\Sigma_{w}}\right)  \geq\text{Length}\left(  P\left(  i\left(
T\right)  \right)  ,g_{\Sigma^{e}}\right)  ,
\]
Where $g_{\Sigma_{w}}$ and $g_{\Sigma^{e}}$ is the induced metrics on
$\Sigma_{w}$ and $\Sigma^{e}.$

Note that $P\left(  i\left(  C_{1}\right)  \right)  =i\left(  C_{1}\right)  ,$
$P\left(  i\left(  C_{3}\right)  \right)  =i\left(  C_{3}\right)  ,$ $P\left(
i\left(  A_{2}\right)  \right)  =i\left(  A_{2}\right)  ,$ there is an
isometry $i_{0}:$ $\Sigma^{e}\longrightarrow\left(  \mathbb{R}^{2}%
,dx^{2}\right)  $ such that $i_{0}\circ P\circ i\left(  C_{1}\right)  =C_{1},$
$i_{0}\circ P\circ i\left(  C_{3}\right)  =C_{3},$ $i_{0}\circ P\circ i\left(
A_{2}\right)  =A_{2}.$ Apply Corollary \ref{Minimal tree}, we have%
\[
\text{Length}\left(  P\left(  i\left(  T\right)  \right)  ,g_{\Sigma
^{e}}\right)  =\text{Length}\left(  i_{0}\circ P\circ i\left(  T\right)
,dx^{2}\right)  >\text{Length}\left(  T,g_{\delta_{n}}\right)  .
\]
Thus we arrive at%
\[
\text{Length}\left(  T,g_{\delta_{n}}\right)  >\text{Length}\left(
T,g_{\delta_{n}}\right)  .
\]
This contradiction finishes the proof of the above proposition.
\end{proof}

Now we give the constructive proof of Theorem 1.1.

\begin{proof}
Step1. Let $\widetilde{k}$ be a smooth function in $\mathbb{R}^{2}$ satisfying%
\begin{align*}
\widetilde{k} &  <0\;\;\;\text{in }B^{n}=B_{2^{-2n}}\left(  2^{-n},0\right)
,\;\;n=1,2,3,\cdots\\
\widetilde{k} &  =0\;\;\;\;\text{otherwise.}%
\end{align*}
Let $u_{1}$ be a smooth solution of%
\[
\triangle u_{1}=-\widetilde{k}.
\]
Then the Gaussian curvature of the metric $g_{1}=e^{2u_{1}}dx^{2}$ satisfies%
\begin{align*}
K_{g_{1}} &  =-e^{-2u_{1}}\triangle u_{1}<0\;\;\;\text{in\ \ \ }B^{n}\\
K_{g_{1}} &  =0\;\;\;\;\;\;\;\;\;\;\;\;\;\;\;\;\;\;\;\text{otherwise.}%
\end{align*}
Step2. Choose a sequence $z_{n,k}$ outside each $B^{n}$ and $\left\{  \left(
x_{1},x_{2}\right)  |x_{2}=0\right\}  $ such that%
\begin{align*}
&  \lim_{k\rightarrow\infty}z_{n,k}\in\partial B^{n}\\
&  \partial B^{n}\subset\overline{\left\{  z_{n,k}\right\}  _{k=1}^{\infty}.}%
\end{align*}
For each $z_{n,k},$ choose a simply connected thin tail $T_{n,k}$ with
$T_{n,k}$ connecting $z_{n,k}$ and the boundary $\partial B_{1}$ such that%
\begin{align*}
&  z_{n,k}\in T_{n,k}\\
&  \partial T_{n,k}\cap\partial B_{1}=\;\;\text{ a piece of arc with positive
length}\\
&  T_{n,k}\subset\mathbb{R}_{+}^{2}=\left\{  \left(  x_{1},x_{2}\right)
|x_{2}>0\right\}  \;\;\;\text{for \ \ \ \ }x_{2}\left(  z_{n,k}\right)  >0\\
&  T_{n,k}\subset\mathbb{R}_{-}^{2}=\left\{  \left(  x_{1},x_{2}\right)
|x_{2}<0\right\}  \;\;\;\text{for \ \ \ \ }x_{2}\left(  z_{n,k}\right)  <0\\
&  T_{n,k}\cap T_{m,j}=\varnothing\;\;\;\text{for }\left(  n,k\right)
\neq\left(  m,j\right)  .
\end{align*}

%
%
%

%

\begin{figure}
[ptbh]
\begin{center}
\includegraphics[
height=7.7431cm,
width=8.1012cm
]%
{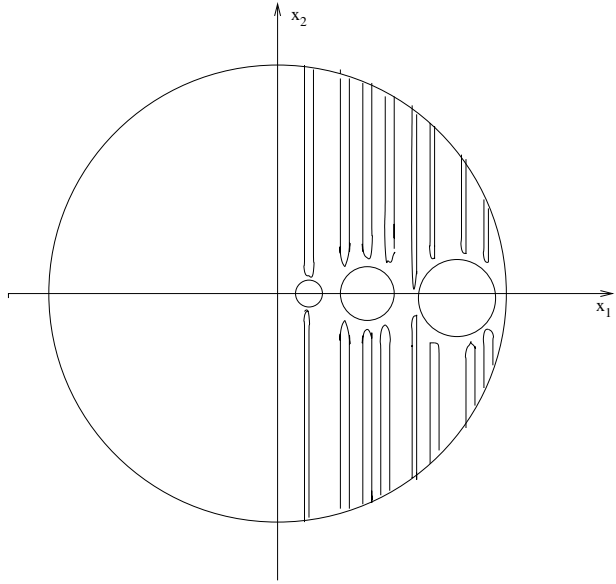}%
\caption{Tails extending to the boundary.}%
\end{center}
\end{figure}

\begin{figure}
[ptbh]
\begin{center}
\includegraphics[
trim=0.000000in 0.000000in -0.001787in 0.004760in,
height=6.0934cm,
width=3.4707cm
]%
{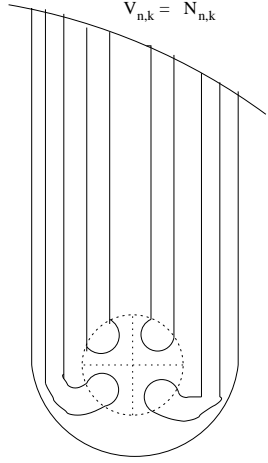}%
\caption{``Details'' of one tail.}%
\end{center}
\end{figure}

We modify the metric $g_{1}$ $=e^{2u_{1}}dx^{2}$ over each tail $T_{n,k}.$ But
we proceed with the tails in the upper and lower half planes separately.

Since $K_{g_{1}}\equiv0$ in the simply connected domain $\mathbb{R}_{+}%
^{2}\backslash\cup_{n=1}^{\infty}B^{n}$. We represent $g_{1}=dy_{+}^{2}$ in
$\mathbb{R}_{+}^{2}\backslash\cup_{n=1}^{\infty}B^{n}$ by a different
coordinate system $y_{+}.$ Over each $T_{n,k}\subset R_{+}^{2},$ we plant a
metric%
\[
g_{2}=e^{2V_{n,k}}dy_{+}^{2}\;\;\;\;\text{in \ \ \ }x^{-1}\left(
T_{n,k}\right)  ,
\]
where $V_{n,k}$ is similar to the one in the construction before Proposition
\ref{II=0}, but the 360 disjoint sub-tails extend to the boundary
$x^{-1}\left(  \partial B_{1}\right)  $ within $x^{-1}\left(  T_{n,k}\right)
.$ We know $V_{n,k}=0$ in $x^{-1}\left(  B_{1}\backslash T_{n,k}\right)  .$
With $V_{n.k}=N_{n,k}$ chosen large enough on $x^{-1}\left(  \partial
B_{1}\right)  $\ intersection with the $x$ pre-image of the 360 sub-tails, we
make%
\[
\triangle V_{n,k}\geq0\;\;\;\text{in\ \ \ }x^{-1}\left(  B_{1}\right)  .
\]

We modify the metric $g_{1}=e^{2u_{1}}dx^{2}$ over the tails in the lower half
plane $\mathbb{R}_{-}^{2}$ with different coordinate system in the same way.

So far, we obtain a new metric $g_{2}=e^{2u_{2}}dx^{2}$\ in $B_{1}$ (which may
not be smooth). We modify $g_{2}$ over the tails one last time.

Let%
\begin{align*}
g_{3} &  =e^{2\epsilon_{n,k}V_{n.k}}dy_{+}^{2}\;\;\;\;\;\;\text{in
\ \ \ \ }x^{-1}\left(  T_{n,k}\right)  \text{\ \ \ for \ \ \ \ }T_{n,k}%
\subset\mathbb{R}_{+}^{2}\\
g_{3} &  =e^{2\epsilon_{n,k}V_{n.k}}dy_{-}^{2}\;\;\;\;\;\;\text{in
\ \ \ \ \ }x^{-1}\left(  T_{n,k}\right)  \text{\ \ \ for \ \ \ \ \ \ }%
T_{n,k}\subset\mathbb{R}_{-}^{2}.
\end{align*}
By choosing $\epsilon_{n,k}>0,$ $\epsilon_{n,k}\longrightarrow0$ sufficiently
fast for $k\longrightarrow\infty$, we can assure $g_{3}=e^{2u_{3}}$ $dx^{2}$
is a smooth metric with $K_{g_{3}}\leq0$\ in $B_{1}.$

Step 3. Suppose there is an isometric embedding%
\[
i:\left(  B_{r},g\right)  \longrightarrow R^{3}%
\]
for some $r>0.$ Then there is $n_{\ast}$ such that%
\[
B^{n_{\ast}}\subset B_{r}.
\]
Applying Proposition \ref{II=0}, we have%
\[
II\circ i=0\;\;\;\;\text{on\ \ \ }\partial B^{n_{\ast}}.
\]
We may assume $i\left(  B_{r}\right)  $ is represented as a graph
$x_{3}=f\left(  x_{1},x_{2}\right)  $ with $\nabla f\left(  0,0\right)  =0.$
Also we may assume the projection of $i\left(  B^{n_{\ast}}\right)  $ down to
x$_{1}$\_x$_{2}$ plane is a domain $\Omega.$ Then%
\begin{align*}
&  \det D^{2}f=K_{g}\left(  1+\left|  \nabla f\right|  ^{2}\right)
^{2}<0\;\;\;\text{in\ \ \ }\Omega\\
&  D^{2}f=0\;\;\;\;\;\;\;\;\;\;\;\;\;\;\;\;\;\;\;\;\;\;\text{on\ \ \ }%
\partial\Omega.
\end{align*}
From $D^{2}f=0$ on $\partial\Omega,$ it follows that $\nabla f=const$. on
$\partial\Omega$ and $f$ coincides with a linear function on $\partial\Omega.$
After subtracting the linear function from $f$, we may further assume $f=0$ on
$\partial\Omega.$ We still have $\det D^{2}f<0$ in $\Omega.$ From the maximum
principle, we see that $f\equiv0$ in $\Omega.$ This contradiction finishes the
proof of Theorem 1.1.
\end{proof}

\section{Metric with negative curvature except for one point}

Relying on the metric constructed in Section 2, we construct a smooth metric
$g$ in $B_{1}$ with negative Gaussian curvature except for one point, namely,
$K_{g}\left(  x\right)  <0$\ for $x\neq0,$ such that the surface $\left(
B_{1},g\right)  $ is not $C^{3,\alpha}$ isometrically embeddable in
$\mathbb{R}^{3}$ even locally near $0.$

For any surface $\left(  \Omega,g\right)  $, we define the $C^{3,\alpha}$
isometric embedding norm by%
\[
\left\|  \left(  \Omega,g\right)  \right\|  _{E}=\inf\left\{  \left\|
II\left(  i\left(  \Omega\right)  \right)  \right\|  _{C^{1,\alpha}%
}|\;C^{3,\alpha}\;\text{isometric\ embedding}\;i:\left(  \Omega,g\right)
\longrightarrow R^{3}\right\}  .
\]
Now we give a constructive proof of Theorem 1.2.

\begin{proof}
Let the annulus $A^{n}=B_{1/n}\backslash B_{1/\left(  n+1\right)  }\subset
R^{2}.$ We construct a metric $g=e^{2u_{0}}dx^{2}$ on $B_{1}$ such that a
non-isometrically embeddable metric $g$ as in Theorem 1.1 is planted (not just
cut and pasted) over each annulus $A^{n}.$

%
%

\begin{figure}
[ptbh]
\begin{center}
\includegraphics[
trim=0.000000in 0.000000in 0.003832in -0.003626in,
height=7.6794cm,
width=8.1012cm
]%
{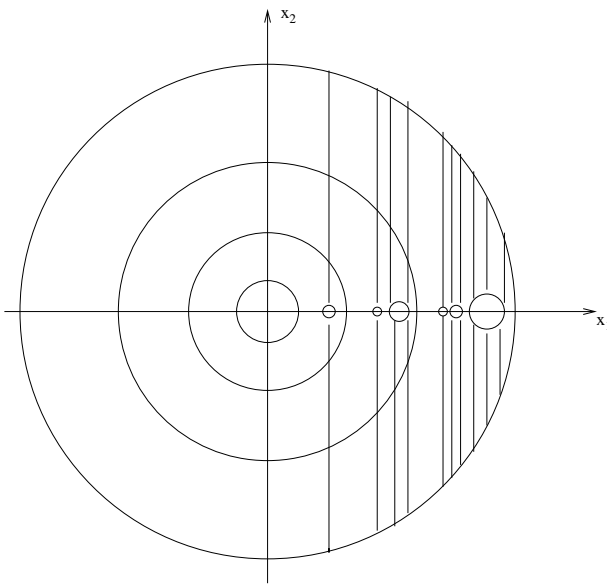}%
\caption{Non-embeddable metric in each annulus.}%
\end{center}
\end{figure}

Set%
\[
\widetilde{\varphi_{n}}\left(  r\right)  =\left\{
\begin{array}
[c]{l}%
e^{-\frac{1}{r-1/n}}\qquad r=\left|  x\right|  >\frac{1}{n}\\
0\qquad\qquad\quad0\leq r\leq\frac{1}{n}%
\end{array}
\right.
\]
We choose $\mu_{1}>0,\;\mu_{2}>0,$ $\cdots,$ $\mu_{n}>0,$ $\cdots$ such that
$\varphi_{n}=\mu_{n}\widetilde{\varphi_{n}}$ satisfies that $\sum
_{n=1}^{\infty}\varphi_{n}$ is smooth and even $\sum_{n=1}^{\infty}%
\epsilon_{n}\varphi_{n}$ is smooth for $\left(  \epsilon_{1},\epsilon
_{2},\cdots\right)  \in l^{\infty}.$

For $\epsilon=\left(  \epsilon_{1},\epsilon_{2},\cdots\right)  \in
l_{+}^{\infty},$ that is $\epsilon_{1}>0,\;\epsilon_{2}>0,\cdots$ and
$\left\|  \epsilon\right\|  _{\infty}=\max\epsilon_{m}<+\infty$ , set%
\begin{align*}
&  \Phi_{\epsilon}=\sum_{m=1}^{\infty}\epsilon_{m}\varphi_{m}\\
&  g_{v}=e^{2\left(  u_{0}+v\right)  }dx^{2}.
\end{align*}
By the construction, $\left(  A^{n},e^{2u_{0}}dx^{2}\right)  $ is not $C^{3}$
isometrically embeddable in $\mathbb{R}^{3}$ for any $n,$ then we have the following.

There exists $0<\eta_{1}$ such that $\left\|  \left(  A^{1},g_{\Phi_{\epsilon
}}\right)  \right\|  _{E}\geq1$ for $\epsilon\in l_{+}^{\infty}$ with
$\left\|  \epsilon\right\|  _{\infty}\leq\eta_{1}.$

Next there exists $0<\eta_{2}<\eta_{1}$ such that $\left\|  \left(
A^{m},g_{\Phi_{\epsilon}}\right)  \right\|  _{E}\geq m$ for $m=1,2$ and
$\epsilon=\left(  \eta_{1},\epsilon_{2},\epsilon_{3},\cdots\right)  \in
l_{+}^{\infty}$ with $\left\|  \left(  0,\epsilon_{2},\epsilon3,\cdots\right)
\right\|  _{\infty}\leq\eta_{2}.$

Inductively there exists $0<\eta_{k}<\eta_{k-1}$ such that $\left\|  \left(
A^{m},g_{\Phi\epsilon}\right)  \right\|  _{E}\geq m$ for $m=1,2,\cdots,k$ and
with $\epsilon=\left(  \eta_{1},\eta_{2},\cdots,\eta_{k},\epsilon
_{k+1},\epsilon_{k+2},\cdots\right)  \in l_{+}^{\infty}$ with $\left\|
\left(  0,\cdots,0,\epsilon_{k+1},\epsilon_{k+2},\cdots\right)  \right\|
_{\infty}\leq\eta_{k}.$

$\cdots$

Finally let $\Psi=\sum_{m=1}^{\infty}\eta_{m}\varphi_{m}$ , $g=g_{\Psi}.$ We
see that
\begin{align*}
&  \left\|  \left(  A^{m},g\right)  \right\|  _{E}\geq m\;\;\;\text{for}%
\;\;m=1,2,3,\cdots\;\;\;\\
&  K_{g}\left(  x\right)  <0\;\;\text{\ for}\;\;x\neq0\;\;\;\;\text{and}\\
&  K_{g}\left(  0\right)  =0.
\end{align*}
It follows that there is no $C^{3,\alpha}$ isometric embedding of $\left(
B_{r}\left(  0\right)  ,g\right)  $ in $\mathbb{R}^{3}$ for any $r>0,\;\alpha>0.$
\end{proof}


\begin{thebibliography}{99}
\bibitem[ C]{C}Cartan, E., \emph{Sur la possibilit\'{e} de plonger un espace
riemannian donn\'{e} dans un espace euclidien,} Ann. Soc. Pol. Math.,
\textbf{5} (1927), 1--7.

\bibitem[ GR] {GR}Gromov, M. L. and Rokhlin, V. A., \emph{Embeddings and
immersions in Riemannian geometry,} Uspehi Mat. Nauk, \textbf{25} (1970)
3--62; English translation in Russian Math. Survey \textbf{25} (1970), 1--57.

\bibitem[ HHL] {HHL}Han, Q., Hong, J. X., and Lin, C. S., to appear.

\bibitem[ HW1] {HW1}Hartman, P. and Winter, P., \emph{On the asymptotic curves
of a surface,} Amer. J. Math., \textbf{73} (1951), 149--172.

\bibitem[ HW2] {HW2}Hartman, P. and Winter, P., \emph{On hyperbolic partial
differential equations,} Amer. J. Math., \textbf{74} (1952), 834--864.

\bibitem[ I] {I}Iwasaki, N., \emph{Applications of Nash-Moser theory to
nonlinear Cauchy problems,} Proc. Sympos. Pure Math., \textbf{45} (1986) Part
1, 525--528.

\bibitem[ H] {H}Hong, J. X., \emph{Cauchy problems for degenerate hyperbolic
Monge-Amp\`{e}re equations and some applications,} J. Partial Differential
Equations \textbf{4} (1991), 1--18.

\bibitem[ J] {J}Janet, M., \emph{Sur la possibilit\'{e} de plonger un espace
riemannian donn\'{e} dans un espace euclidien,} Ann. Soc. Pol. Math.,
\textbf{5} (1926), 38--43.

\bibitem[ K] {K}Kuiper, N. H., \emph{On }$C^{1}$\emph{-isometric embeddings,
I, II,} Indag. Math., \textbf{17} (1955), 545--556, 683--689.

\bibitem[ L1] {L1}Lin, C. S., \emph{The local isometric embedding in
}$\mathbb{R}^{3}$\emph{ of 2-dimensional Riemannian manifolds with nonnegative
curvature, }J. Diff. Goem., \textbf{21} (1985), 213--230.

\bibitem[ L2] {L2}Lin, C. S., \emph{The local isometric embedding in
}$\mathbb{R}^{3}$\emph{ of two-dimensional Riemannian manifolds with Gaussian
Curvature changing sign cleanly,} Comm. Pure Appl. Math., \textbf{39} (1986), 867--887.

\bibitem[ N] {N}Nadirashvili, N., \emph{The local embedding problem for
surfaces,} preprint, 2002.

\bibitem[ Na1] {Na1}Nash, J., $C^{1}$\emph{ isometric imbeddings,} Ann. of
Math., \textbf{60} (1954), 383--396.

\bibitem[ Na2] {Na2}Nash, J., \emph{The imbedding problem for Riemannian
manifolds,} Ann. of Math., \textbf{63} (1956), 20--63.

\bibitem[ Ni] {Ni}Nirenberg, L., \emph{The Weyl and Minkowski problems in
differential geometry in the large,} Comm. Pure Appl. Math., \textbf{6}
(1953), 337--394.

\bibitem[ P1] {P1}Pogorelov, A. V., \emph{Intrinsic estimates for the
derivatives of the radius vector of a point on regular convex surface,} Dokl.
Akad. Nauk SSSR (N.S.), \textbf{66} (1949), 805--808.

\bibitem[ P2] {P2}Pogorelov, A. V., \emph{An example of a two-dimensional
Riemannian metric that does not admit a local realization in E}$^{3},$ Dokl.
Akad. Nauk SSSR, \textbf{198} (1971) 42--43; English translation in Soviet
Math. Dokl. \textbf{12} (1971), 729--730.

\bibitem[ Po] {Po}Poznyak, \`{E}. G., \emph{Isometric embeddings of
two-dimensional Riemannian metrics in Euclidean space,} Uspekhi Mat. Nauk,
\textbf{28} (1973), 47--76; English translation in Russian Math. Survey,
\textbf{28} (1973), 47--77.

\bibitem[ Y] {Y}Yau, S.-Y., \emph{Seminar on Differential Geometry,} Annals of
Math. Studies, \textbf{102}, Princeton, 1982.
\end{thebibliography}
\end{document}